\documentclass[a4paper, 12pt]{scrartcl}

\usepackage[T1]{fontenc}
\usepackage[utf8]{inputenc}
\usepackage{graphicx}
\usepackage[english]{babel}
\usepackage{amsmath}
\usepackage{amsthm}
\usepackage{amssymb}
\usepackage[super]{cite}

\usepackage{algorithm}
\usepackage{algorithmic}

\usepackage{chemfig,chemmacros,chemformula}

\graphicspath{{figures/}}

\newcommand{\R}{\ensuremath{\mathbb R}}
\newcommand{\Rnn}{\ensuremath{\mathbb{R}_{+}}}

\DeclareMathOperator*{\argmin}{arg\,min}

\newcommand{\diag}[1]{\ensuremath{\text{diag}(#1)}}
\newcommand{\abs}[1]{\ensuremath{\text{abs}(#1)}}

\newcommand{\defby}{\mathrel{\mathop:}=}

\newcommand{\bydef}{=\mathrel{\mathop:}}

\newcommand{\printurl}[1]{\url{#1}}

\newcommand{\norm}[1]{\lVert #1 \lVert}

\newcommand{\noiselevel}{\delta}

\theoremstyle{plain}

\theoremstyle{remark}

\theoremstyle{definition}

\makeatletter
\renewcommand\@biblabel[1]{#1.}
\makeatother

\title{Using separable non-negative matrix factorization techniques
for the analysis of time-resolved Raman spectra}

\author{
Robert Luce\thanks{
EPF Lausanne, SB MATHICSE ANCHP, MA B2 525, Station 8, CH-1015
Lausanne, Switzerland, \ttfamily{robert.luce@epfl.ch}
},
Peter Hildebrandt\thanks{
TU Berlin, PC 14, Stra{\ss}e des 17. Juni 135, 10623 Berlin, Germany,
\ttfamily{hildebrandt@chem.tu-berlin.de}
},
Uwe Kuhlmann\thanks{
TU Berlin, PC 14, Stra{\ss}e des 17. Juni 135, 10623 Berlin, Germany,
\ttfamily{uwe.kuhlmann@tu-berlin.de}
},
J\"org Liesen\thanks{
TU Berlin, MA 4-5, Stra{\ss}e des 17. Juni 136, 10623 Berlin, Germany,
\ttfamily{liesen@math.tu-berlin.de}
}
} % End of \author

\begin{document}

\maketitle

\begin{abstract}

The key challenge of time-resolved Raman spectroscopy is the
identification of the constituent species and the analysis of the
kinetics of the underlying reaction network. In this work we present
an integral approach that allows for determining both the component
spectra and the rate constants simultaneously from a series of
vibrational spectra. It is based on an algorithm for non-negative
matrix factorization which is applied to the experimental data set
following a few pre-processing steps. As a prerequisite for physically
unambiguous solutions, each component spectrum must include one
vibrational band that does not significantly interfere with
vibrational bands of other species. The approach is applied to
synthetic ``experimental'' spectra derived from model systems
comprising a set of species with component spectra differing with
respect to their degree of spectral interferences and signal-to-noise
ratios. In each case, the species involved are connected via
monomolecular reaction pathways. The potential and limitations of the
approach for recovering the respective rate constants and component
spectra are discussed.

\end{abstract}

\section{Introduction}

Raman spectroscopy is a versatile tool to probe molecular structure
changes that are associated with the temporal evolution of chemical or
physical processes
\cite{SahooUmapathyParker2011,BalakrishnanEtAl2008,OndriasSimpsonLarson1996}.
Since time-resolved Raman spectroscopy is applicable in a wide dynamic
range down to the femtosecond time scale, it is capable to monitor
quite different events, including intramolecular rearrangements in the
excited state as well as chemical reactions in the ground state.  The
individual Raman spectra measured as a function of time represent a
superposition of the intrinsic spectra of the individual species
or molecular states that are involved in the reaction sequence. The
relative contributions of the various spectra to each measured spectrum
then reflect the actual composition of the sample at the respective
time, and hence the entire series of experimental spectra represents
the kinetics of the underlying processes. While just this
information can also be provided by other transient optical
techniques, the unique advantage of time-resolved Raman spectroscopy
resides in the fact that the molecular structure of the individual
species is encoded in the respective component spectra.  This
allows for identifying intermediate states and characterizing their
structural and electronic properties as a prerequisite for determining
the molecular reaction mechanism.  The central task is, therefore, to
disentangle the series of time-resolved Raman spectra in terms of the
individual component spectra and their temporal evolution.

In the past decades, different concepts were designed for analyzing
sets of complex Raman
spectra~\cite{DoepnerEtAl1996,HendlerShrager1994,HenryHofrichter1992,ShinzawaEtAl2009}.
In most cases, these efforts were spurred by the concomitant vivid
development in Raman and IR spectroscopic investigation of cellular systems
for biological and medical
applications.\cite{TalariEtAl2015,BeljebbarEtAl2009}
These studies are dedicated to classify and to distinguish
microorganisms or to identify pathological tissues.  To achieve these
objectives, it is necessary to determine spectral signatures of a
complex ensemble of biomolecules in a certain environment that are
characteristic of a specific type or state of the cellular system. In
these cases, analytical methods for pattern recognition are required
which, in view of the large number of experimental spectra, are based
on statistical procedures, such as principal component analysis,
factor analysis, or singular value
decomposition\cite{DoepnerEtAl1996,HendlerShrager1994,HenryHofrichter1992,ShinzawaEtAl2009,TalariEtAl2015,BeljebbarEtAl2009,WeakleyEtAl2012,ZhangEtAl2005}.
Such approaches provide mathematical solutions but typically not the
intrinsic spectra of the large number of pure components. Thus,
extending multivariate analyses to a series of spectra reflecting
physical or chemical changes of well-defined molecular species may
just afford the number of the species involved but not necessarily
their component spectra.

Such systems, on the other hand, are frequently treated by
least-square methods in which either single Lorentzian/Gaussian bands
and complete component spectra (component analysis) are employed to
achieve a global fit to all experimental
spectra\cite{DoepnerEtAl1996}. The component analysis is quite robust
as the number of degrees of freedom in the fitting process just
corresponds to the number of components, given that all component
spectra are known a priori. If, however, this is not the case, the
“intuition” factor gains weight and thus the overall error increases
substantially with the number of unknown component spectra.

In this work, we tried to overcome these drawbacks by developing an
unsupervised analytical method that is based on non-negative matrix
factorization (NMF).   Such NMF techniques have rarely been used for
factor analysis, with two notable
exceptions\cite{AndoHamaguchi2013,Sawall2010}.  Unlike these previous
approaches, our method takes advantage of the so-called separability
inherent to the measurement data.

This paper is organized as follows.  In Section~\ref{sec:background}
we present the mathematical model and describe the so-called
separability condition that is derived from the specific properties of
complex Raman spectra composed by a finite number of component
spectra. Section~\ref{sec:overall_algo} then describes the numerical
method for computing the NMF and extracting the reaction coefficients.
In Section~\ref{sec:study} we illustrate on a sequence of artificial
first-order reactions the reliability of our method under interference
among the component spectra and under measurement noise. Concluding
remarks are given in Section~\ref{sec:concl}.

\section{Mathematical background}
\label{sec:background}

\subsection{Model for time-resolved Raman spectra of chemical
reactions}
\label{sec:model}

Mathematically, the acquired measurements correspond to certain convex
combinations of sums of Lorentz functions or Lorentzians
that constitute the Raman
spectra of the individual reactants. We will formalize our notion of
the model in this section.

A Lorentzian $L_{x_0, \gamma, I}(x)$ is a non-negative ``peak function''
with its maximum at the base point $x_0\in\R$ (corresponding to the frequency 
of the normal mode), the width at half-height $\gamma>0$ and the 
intensity $I>0$, which is defined by
\begin{equation}
\label{eqn:Lorentz_def}
L_{x_0, \gamma, I}(x) = I \tfrac{\gamma^2}{ (x-x_0)^2 + \gamma^2 }.
\end{equation}
For simplicity we will usually skip the parameters $x_0, \gamma, I$ 
and simply write $L(x)$.

Consider a chemical reaction with $r$ reactant species. Then the Raman spectrum 
$w_s$ of each reactant can be modeled as a non-negative sum of $q_s$ Lorentzians, 
so that
\begin{equation}
\label{eqn:w}
w_s(x) = \sum_{k=1}^{q_s} L^s_k (x),\quad s=1,\dots,r.
\end{equation}
We assume that all base points (or ``peaks'') are located in the finite interval 
$[f_l,f_u]\subset \Rnn\defby [0,\infty)$. Note that one could also implement 
a model where the individual spectra are given by sums of Lorentzians 
and Gaussians, or Gaussians only.

Now we will consider the relative concentrations of the reactant
species.  We denote the concentration function of species $s$ by
\begin{equation*}
	h_s : [0, T] \rightarrow [0,1],\quad s=1,\dots,r,
\end{equation*}
so that $h_s(t)$ corresponds to the relative concentration of species
$s$ at time $t\in [0,T]$ of the reaction.  Consequently, at each time $t$ the
$r$ concentrations sum to $1.0$.  The functions $h_s(t)$ represent the
\emph{reaction kinetics}.

In this functional setting, the time-resolved vibrational Raman
spectrum of the reaction can be modeled as
\begin{equation}
\label{eqn:measurements_functional}
M(x, t)
= \sum_{s=1}^r w_s(x) h_s(t)
= \sum_{s=1}^r \left( \sum_{k=1}^{q_s} L^s_k (x) \right) h_s(t).
\end{equation}
Discretizing \eqref{eqn:measurements_functional} over time through $0
= t_0 < \cdots < t_{n-1} = T$ in $n$ time steps, and over frequencies
through  $f_l = x_1 < \cdots < x_{m} = f_u$ in
$m$ frequencies, models the measured data from the reaction.  We denote
the resulting \emph{measurement matrix} by $M = [m_{ij}] =
[M(x_i, t_{j+1})]\in\Rnn^{m,n}$.

Now we discretize the functions $w_s(x)$ over the frequencies and 
obtain the vectors 
$$w_s=[w_s(x_1),\dots,w_s(x_m)]^T\in\Rnn^m,\quad  s=1,\dots,r,$$
as well as the matrix $W=[w_1,\dots,w_r]\in\Rnn^{m,r}$. Similarly, 
we discretize the functions $h_s(t)$ over time and obtain the vectors 
$$h_s=[h_s(t_0),\dots,h_s(t_{n-1})]^T\in\Rnn^n,\quad  s=1,\dots,r,$$ 
as well as the matrix
$$H=\begin{bmatrix}h_1^T \\ \vdots \\h_r^T\end{bmatrix}\in\Rnn^{r,n}.$$
Then the measurement matrix $M$, which includes all
experimental spectra, can be written as
\begin{equation}
\label{eqn:reaction_mwh}
M = \sum_{s=1}^r w_s h_s^T = WH,
\end{equation}
i.e., the entry-wise non-negative matrix $M$ is the product of the two
entry-wise non-negative matrices $W$ and $H$.

\subsection{NMF for Raman experimental data}
\label{sec:nmf_for_raman}

The task of analyzing time-resolved Raman spectroscopic data
consists of (1) determining the spectra of the individual species
(component spectra), and (2) identifying the underlying reaction kinetics.
Using the notation from Section~\ref{sec:model}, the corresponding
mathematical problem is as follows: Given the non-negative
data matrix $M\in \Rnn^{m,n}$, and assuming that the reaction involves
$r$ species, find non-negative matrices $W \in \Rnn^{m,r}$ and $H \in
\Rnn^{r,n}$ such that
\begin{equation}
\label{eqn:exact_nmf}
    M = WH.
\end{equation}
Note that a factorization of the form \eqref{eqn:exact_nmf} where $W$
and/or $H$ have negative entries has no physical meaning, as neither 
a measured intensity nor a relative concentration can be negative.

The mathematical task of finding a \emph{non-negative matrix
factorization (NMF)} of $M$ is one of formidable difficulty:  Without
further assumptions on the given data, the problem is ill-posed and
its solutions are non-unique in general.  From a complexity point of
view, NMF is NP-hard\cite{Vavasis:2008}. Moreover, an \emph{exact}
factorization as in \eqref{eqn:exact_nmf} is more a theoretical desire
than achievable in practice.  The presence of noise or other forms of
data uncertainty may simply rule out the existence of such a
factorization.

The theory and computation of non-negative matrix factorizations is a
very active research topic, whose popularity gained much from an
article by Lee and Seung on the use of NMF for feature
extraction\cite{LeeSeung:1999}; an earlier work on NMF dates back to
1974\cite{Thomas:1974}.  A useful overview is given in the recent
survey by Gillis\cite{Gillis:HowAndWhy}.

A practically much better justified view is adopted by considering NMF
as an \emph{approximation} problem rather than an exact factorization.
The usual way to give a formal definition of this approximation
problem is as follows:  Let $\norm{\cdot}$ be a matrix norm. Given
the non-negative data matrix $M$, we seek non-negative matrices $W$ and $H$ such
that
\begin{equation*}
\norm{M - WH}
\end{equation*}
is small.  The norm plays the role of a ``distance function'' and
measures how close the eventually found factors $W$, $H$ reproduce the
given data. In this work we mostly use the Frobenius norm, which for
any (rectangular) matrix $A=[a_{ij}]$ is defined by
$\|A\|_F=\sqrt{\sum_{i,j} |a_{ij}|^2}$.

When treating time-resolved Raman data on chemical reactions, another
problem of interest arises.  Usually one does not know \emph{the
number} of intermediate species in the reaction.  Thus the
required information to set up an NMF problem, $r$, is missing.
To overcome this shortcoming, one may assume that $r$ is not very
large (in practice it is often between two and ten).  It then is a
simple matter of trying different values for $r$ and selecting the
best solution (see Section~\ref{sec:err_comp}).  Another heuristic
estimate is available by the number of ``large'' singular values of
the matrix $M$\cite{MalBook02}.

\subsection{The separability condition}
\label{sec:sepnmf}

While the general NMF problem introduced in
Section~\ref{sec:nmf_for_raman} is very difficult in general, there is
a very important special case, the so-called \emph{separable NMF}
problem~\cite{DonohoStodden:2003, AroraEtAl:2012}.  As before, the
measurement matrix is denoted by $M \in \Rnn^{m,n}$, and the entries
of each row of $M$ sum to $1.0$. (This can always be achieved by
applying a diagonal scaling matrix from the left.)  Algebraically, the
data matrix $M$ is called $r$-\emph{separable}, if it can be written
in the form
\begin{equation}
\label{eqn:sepnmf}
    M = W H = Q
    \begin{bmatrix}
        I_r\\
        W'
    \end{bmatrix}
    H,
\end{equation}
where $I_r$ is the $r$-by-$r$ identity matrix, and $Q \in \R^{m,m}$ is a
permutation matrix.  Separability implies that all the rows of $M$ can
be reconstructed by using only $r$ rows of $M$ (these constitute the
factor $H$) by convex combinations with weights given through $W'$.

The separability condition in equation~\eqref{eqn:sepnmf} can be
interpreted in the model for time-resolved Raman spectroscopy data
from the previous section as follows. The measurement matrix $M$ is
separable when each species $s$, represented by the $s$-th column of
$W$, has a \emph{characteristic frequency} $x_s$ at which $w_s(x_s) >
0$ (see~\eqref{eqn:w}), but $w_{\tilde{s}}(x_s) = 0$ for
all other species $\tilde{s} \neq s$.
If such a frequency is present for each species, the
measurement matrix $M$ contains rows that are equal to the rows of the
sought-for kinetic matrix $H$.  Thus, the primary task is to determine
the characteristic frequencies.  For example, if all the species
involved in the reaction contain a Lorentz band that does not
interfere with a Lorentz band of any other species, then the
separability condition is satisfied.

However, recall from the definition \eqref{eqn:Lorentz_def} that the
intensity of a Lorentzian $L_{x_0, \gamma, I}(x)$ vanishes reciprocally
to a quadratic function in the distance from the base frequency $x_0$.
In particular, the separability condition as explained above
\emph{cannot} be satisfied in any case, since $w_s(x) > 0$ for every
species $s$ and frequency $x$.  Consequently, the
factorization~\eqref{eqn:reaction_mwh} will not exactly be
separable (in the theoretical framework of Section~\ref{sec:model}).

While an \emph{exact} interference-free set of characteristic
frequencies is in general impossible, the interference can be
\emph{numerically} small or even zero if the corresponding base
frequencies of the characteristic bands are sufficiently far apart.
Algebraically, it means that instead of the exact separability
condition~\eqref{eqn:sepnmf}, our problem is properly described by a
\emph{near-separable} problem, meaning that
\begin{equation}
\label{eqn:near_sep_interfere}
M = WH = Q
\begin{bmatrix}
I_r + R\\
W'
\end{bmatrix}
H
=
Q
\begin{bmatrix}
I_r\\
W'
\end{bmatrix}
H
+
\underbrace{
Q
\begin{bmatrix}
R\\
0
\end{bmatrix}
H
}_{\bydef N_I}
=
Q
\begin{bmatrix}
I_r\\
W'
\end{bmatrix}
H
+
N_I.
\end{equation}
Here we interpret the matrix $N_I$ as \emph{noise} originating from
interference of the species at frequencies at which some species is
strongly dominant in comparison with the other species.  If $\|N_I\|$ is
small, the original factors $W, H$ may be well approximated by
applying an algorithm for separable NMF to $M$.
Quantitative investigations on the allowable noise for some algorithms
have been persued in a purely mathematical
context.\cite{GillisVavasis:2012, Gillis:2013, GillisLuce:2013}

In Section~\ref{sec:move_bands}, we provide a
numerical study on a model problem that shows the effect of growing
$\norm{N_I}$ on the overall approximation quality.

So far we considered only ideal, noise free data measurements.  
Real experimental data involve measurement noise, and hence
\begin{equation}
\label{eqn:add_noise}
M = WH + N_M,
\end{equation}
where we assume the noise to be purely additive, i.e. $N_M$ is
componentwise non-negative.
Algebraically, the measurement noise is no different from the noise
arising from interference.  The effect of measurement noise is studied
numerically in Section~\ref{sec:noisebench}.

Note that the separability assumption is widely used in other fields,
such as hyperspectral imaging, text mining, or other blind source
separation
applications~\cite{KumarSindhwaniPrabhanjan2013,Bioucas-DiasEtAl2012,ChanEtAl2008}.
In some of these contexts, the separability assumption is assumed to
be satisfied with respect to the time axis (in our notation, for
$M^T$).  In our application this would mean that for each species
there exists some point in time at which the relative concentration of
the species is $1.0$, which is highly unlikely for a typical reaction.

\section{Computational method}

\label{sec:overall_algo}
\label{sec:comp_aspects}

The method we use for the identification of the reaction kinetics is
based on the successive non-negative projection algorithm
(SNPA)~\cite{Gillis2014}.  SNPA is an algorithm for computing the factor $H$
in~\eqref{eqn:sepnmf} (the kinetics), provided that the problem at
hand is near-separable.  In the language of Section~\ref{sec:sepnmf},
we use SNPA to compute approximate characteristic frequencies of the
species.  We have chosen SNPA for its computational
speed and robustness with respect to noise, but any other algorithm
for separable NMF could be used as well.

While SNPA is at the center of our method, we also need to deal with a
number of other computational tasks.  The overall method is shown in
Algorithm~1.

\begin{algorithm}[h]
\caption{Species and kinetics identification via separable NMF}
\begin{algorithmic}[1]
\REQUIRE Data matrix $M$, number of species $r$
\ENSURE Approximate species $W$ and kinetics $H$ s.t. $M \approx WH$;
reaction coefficients $K$.
\medskip
\STATE Filter out noisy rows (frequencies) of $M$.
\STATE Smooth measurements in direction of the columns of $M$ (time).
\STATE Apply SNPA to $M$ to obtain pseudo-kinetic $\hat{H}$.
\STATE Scale $H \leftarrow D
\hat{H}$ so that the columns of $H$ sum approximately to 1.0.
\STATE Compute corresponding spectra $W$ such that $M \approx
WH$.
\STATE Extract reaction coefficient matrix $K$ from kinetic $H$.
\end{algorithmic}
\end{algorithm}

\paragraph*{Step 1: Removing insignificant frequencies}

In our approach it has turned out to be useful to remove all the
experimental data (intensity-frequency pairs) at frequencies which did
not display any significant intensities with respect to the
measurement noise level.  Here, we first estimate the noise level
inherent to the measurement data by the standard deviation on the
frequency having the least mean intensity.  If we make the reasonable
assumption that there are frequencies at which \emph{only} noise is
measured, such a frequency will have minimum mean intensity,
and its standard deviation is a measure for the noise level.
Algebraically, this filtering just leads to removal of some rows of
$M$.  In order to simplify the notation, we will still assume that $M
\in \Rnn^{m,n}$.

\paragraph*{Step 2: Smoothing the data} An useful preprocessing step
is to ``smooth'' the measurement data in direction of the time.  In all
numerical experiments in the following section that involve
measurement noise, we smoothed the data using a running mean with a 
window size of $5$.  Algebraically, this smoothing just effects that
each data entry $m_{ij}$ of $M$ is replaced by the mean of the values
$m_{ik}$ for $j-2 \le k \le j+2$.

\paragraph*{Step 3: Find characteristic frequencies}

Subsequently, we use SNPA to find a set of $r$ approximate
characteristic frequencies.  If $\mathcal{K} \subset
\{1,\dotsc,m\}$ is the set of $r$ indices computed by SNPA, we
obtain the pseudo-kinetic matrix $\hat{H} \in \Rnn^{r,n}$ by 
stacking the $r$ rows of $M$ indexed by $\mathcal{K}$ (i.e., 
$\hat{H} = M(\mathcal{K}, :)$ in Matlab-like notation).  Note that the columns of
$\hat{H}$ may not sum to (approximately) $1.0$, and hence $\hat{H}$
cannot be interpreted as a reaction kinetic matrix.

\paragraph*{Step 4: Scaling $\hat{H}$}

In order to obtain a reaction kinetic matrix whose columns sum
approximately to $1.0$, we next compute a diagonal scaling matrix $D
\in \Rnn^{r,r}$ such that $D\hat{H}$ has this property.  Here use 
that the approximation error of an NMF is invariant under such 
scalings, viz. $WH = (WD^{-1})(DH)$.
To find a suitable scaling matrix $D$, we solve the
non-negative least squares problem
\begin{equation*}
    \min_{d \in \Rnn^r} \norm{\hat{H}^T d - e}_F,
\end{equation*}
where $e =[1,\dotsc,1]^T \in \R^n$.  After finding
the optimal scaling values $d_1, \dotsc, d_r$, we
rescale the kinetic matrix $H = \diag{d}\hat{H}$.  In our experiments
described in the following section, we used Matlab's {\ttfamily{lsqnonneg}}
function.

\paragraph*{Step 5: Compute corresponding spectra $W$}

With the kinetic matrix $H$ and the measurement data
matrix $M$, we now determine the spectra (i.e., the factor $W$
in~\eqref{eqn:reaction_mwh}), which can be computed by
solving the convex minimization problem
\begin{equation*}
	\min_{W \in \Rnn^{m,r}} \; \norm{M - WH}_F^2
\end{equation*}
using standard techniques, e.g.~\cite{VanBenthemKeenan:2004}.

\paragraph*{Step 6: Extracting reaction coefficients}

To extract the reaction coefficients (rate constants) from a given
kinetic matrix $H$, we restrict the analysis to the case that all
reaction steps are of first order.  Recall from
Section~\ref{sec:model} that if $H$ is the reaction kinetics matrix of
the true kinetics function $h(t) = [h_1(t), \dotsc, h_r(t)]^T$, and $K
\in \R^{r,r}$ is the matrix of reaction coefficients, then $h(t) =
e^{K t} h_0$, where $h_0$ is the initial concentration vector.  If
$K$ is not known, it can be recovered from $H$ by solving the
nonlinear least squares problem
\begin{equation*}
    \argmin_{K \in R^{r,r}} \sum_{j=1}^n \norm{H_j - e^{Kt_j} h_0}_F,
\end{equation*}
where $H_j \in \R^{r}$ denotes the $j$-th column of $H$.  If $H$ has 
full row rank, then $K$ is uniquely determined by $H$.  In our experiments 
we used Matlab's {\ttfamily{fminunc}} function.

\section{Numerical study with synthetic data}
\label{sec:study}

In this section we illustrate the effectiveness of our method using a
sequence of artificial first-order reactions involving five species.
We first describe the model reactions and the component spectra of the
involved species (``fingerprints'') in
Section~\ref{sec:synthetic_reaction}.  The results in
Section~\ref{sec:noiseless} show that in the low-interference,
noiseless regime, both the kinetics and species fingerprints are
perfectly recovered.  Then, in Section~\ref{sec:move_bands}, we study
numerically the effect of increasing interference among the species,
and in Section~\ref{sec:noisebench} we also add measurement
noise to the data and study the recovery quality.  Finally in
Section~\ref{sec:err_comp}, we address the question of determining
the correct number of species.

\subsection{Description of the model reaction}
\label{sec:synthetic_reaction}

The reaction scheme is set-up by five species A, B, C, D and E which
are inter-related by first-order reactions.  These reactions are
characterized by reaction coefficients (in arbitrary units of
reciprocal time) as given as follows:

\begin{center}
\schemestart
A
\arrow(--bb){<=>[0.53][0.02]}
B
\arrow{<=>[0.43][0.25]}
C
\arrow{<=>[0.11][0.1]}
E
\arrow(@bb--){->[*{0}0.21]}[-90]
D
\schemestop
\end{center}

\begin{figure}[t]
    \begin{center}
        \includegraphics[width=.48\textwidth]{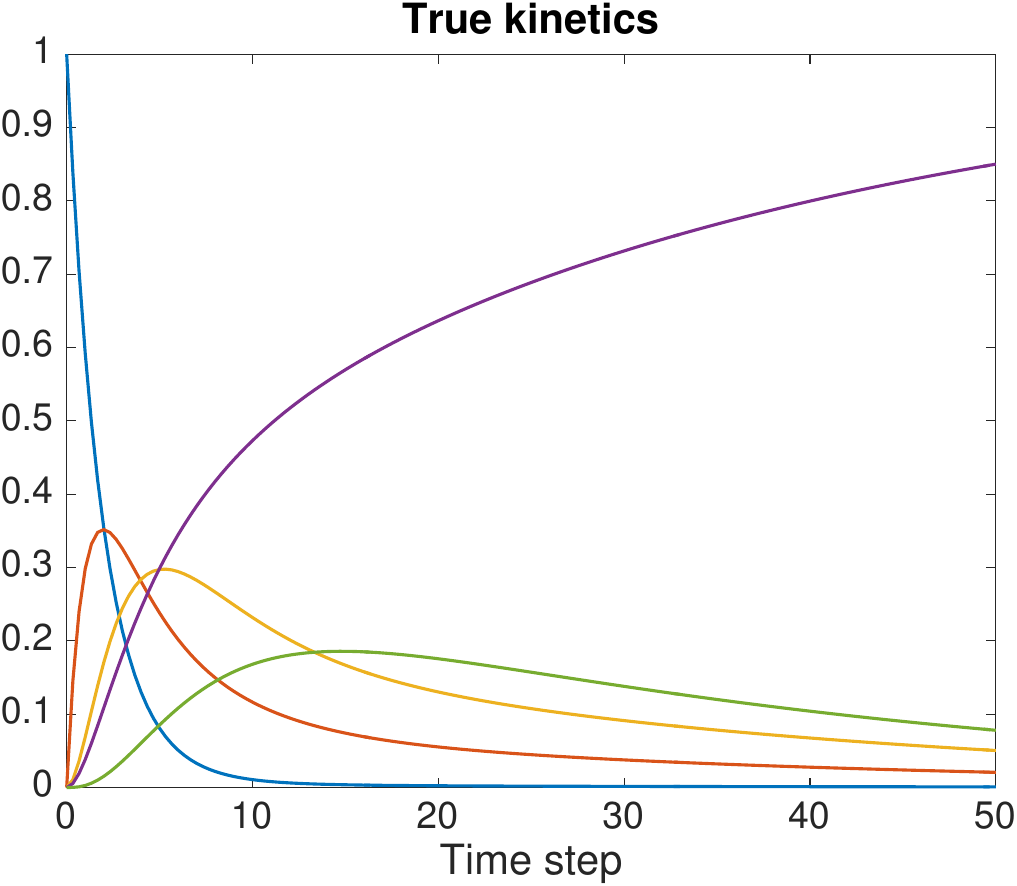}
        \hfill
        \includegraphics[width=.48\textwidth]{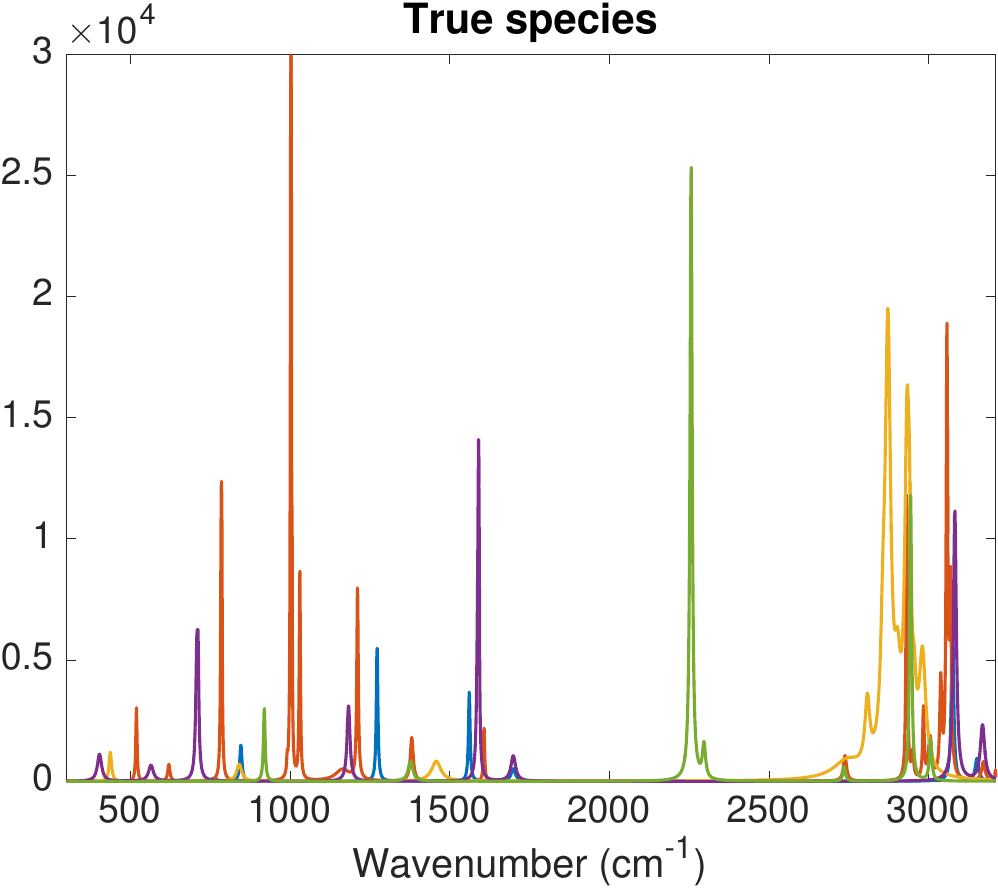}
    \end{center}
    \caption{Artificial sequence of first order reactions with five
species.  The kinetics are shown in the left panel and the species
fingerprints (component spectra) are displayed in the right panel.
The resulting measurement data are shown in
Figure~\ref{fig:synreaction_datavis} (top).\label{fig:synreaction}}
\end{figure}

We let species A be the only educt in the reaction, resulting in the
initial concentration vector $h_0 \defby h(t_0) = [1, 0, 0, 0, 0]^T$.
If we denote the reaction coefficient matrix corresponding to the
above reaction scheme by
\begin{equation}
\label{eqn:k_matrix}
K =
\begin{bmatrix}
    -0.53 & 0.02 &  0.0  &  0 &  0.0\\
     0.53 &-0.66 &  0.25 &  0 &  0.0\\
     0.0  & 0.43 & -0.36 &  0 &  0.1\\
     0.0  & 0.21 &  0.0  &  0 &  0.0\\
     0.0  & 0.0  &  0.11 &  0 & -0.1
\end{bmatrix},
\end{equation}
the reaction kinetics are given as a function
over time by 
$$h(t) = [h_1(t),\dots,h_5(t)]^T=e^{Kt} h_0,$$
as displayed in Figure \ref{fig:synreaction} (left).  
The corresponding kinetics matrix
$H$ is obtained by discretization of $h(t)$ at equidistant steps
$t_0, \dotsc, t_{n-1}$, so that $H = [h(t_0), \dotsc, h(t_{n-1})]$.

The five component spectra are constructed by arbitrarily chosen
combinations of Lorentzians, inspired by the Raman spectra of
various organic compounds.  These component spectra were constructed
such that each species has at least one characteristic frequency, so
that the imposed separability condition (see Section~\ref{sec:sepnmf})
is satisfied.  The fingerprints constitute the columns of the matrix
$W$ (Figure~\ref{fig:synreaction} (right)).  Visual inspection already
reveals that each of the five species has at least one characteristic
frequency.

\begin{figure}[p]
    \begin{center}
        \includegraphics[width=.80\textwidth]{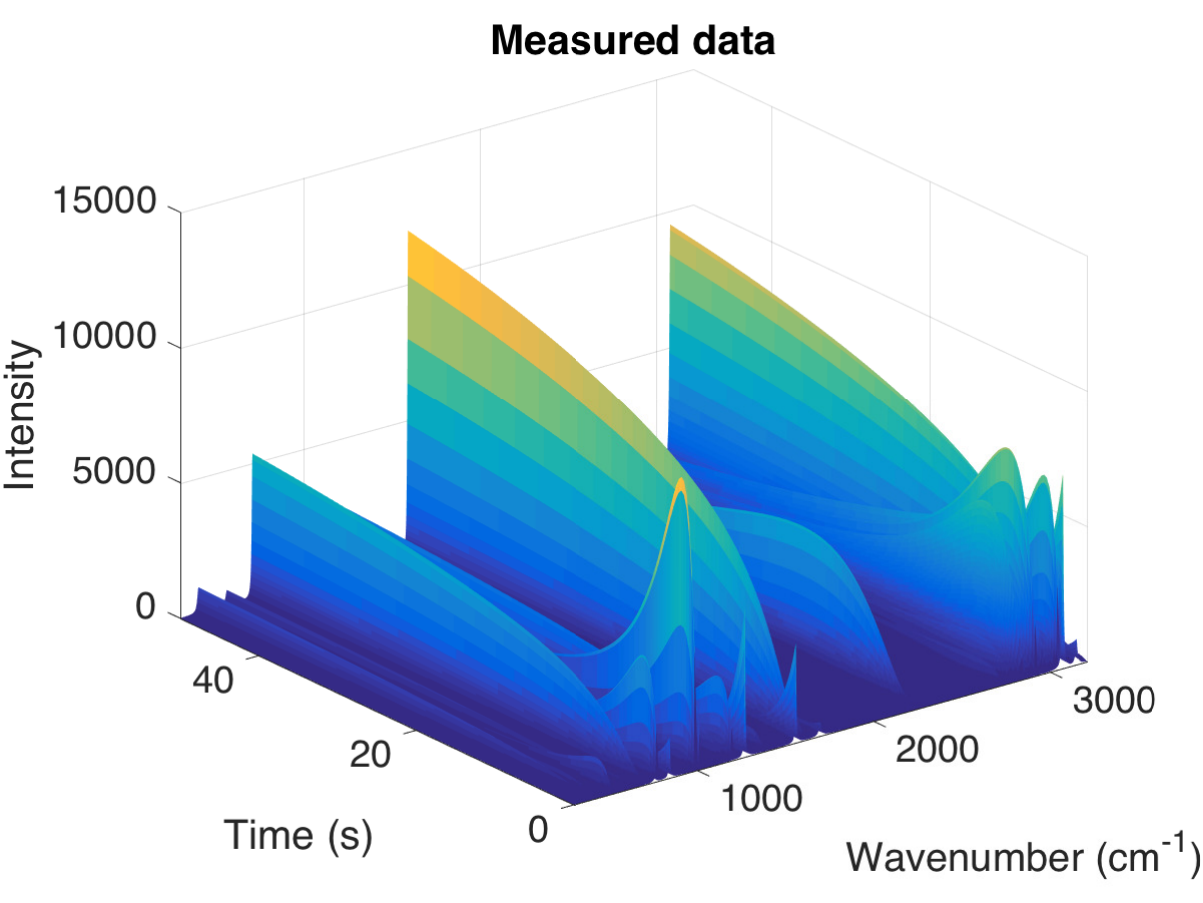}\\
        \includegraphics[width=.80\textwidth]{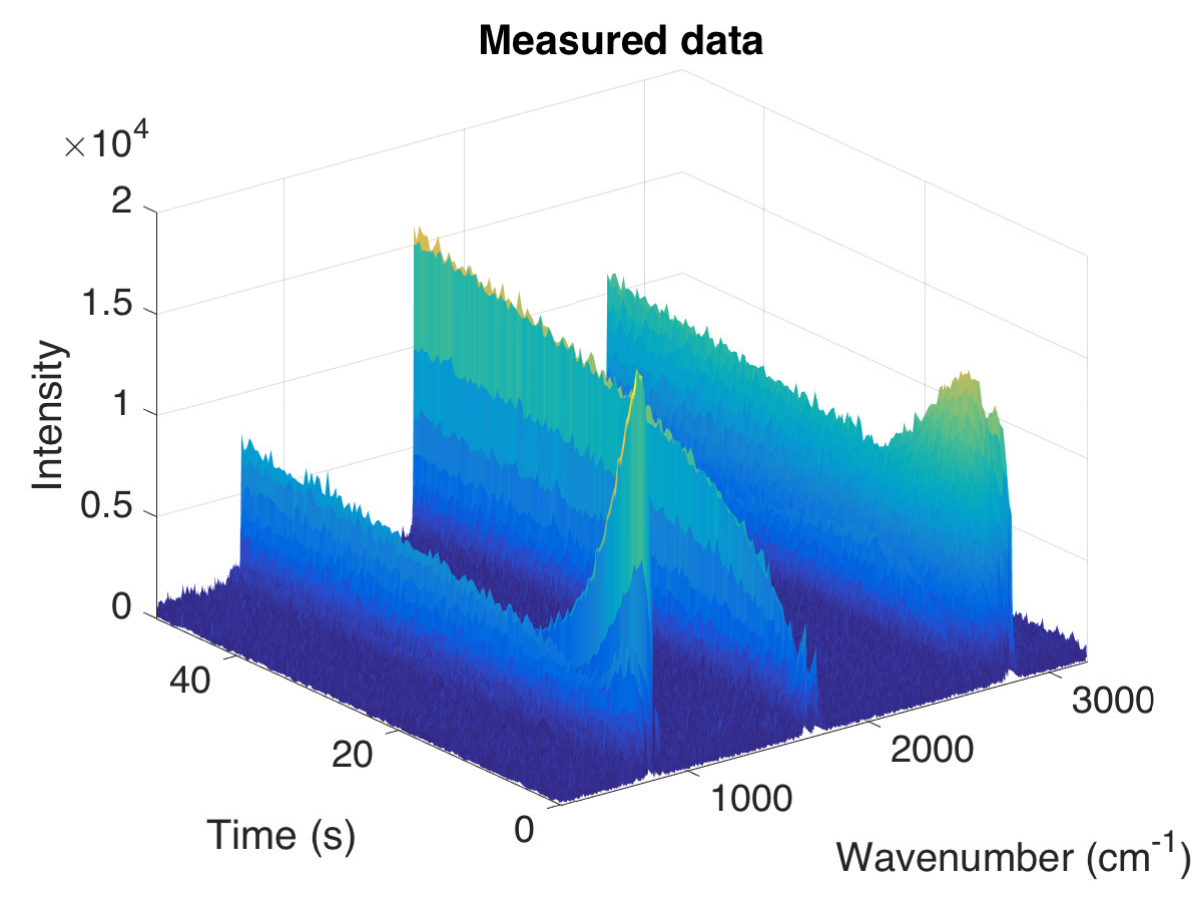}
    \end{center}
    \caption{Visualization of the measurement data matrix $M$ for the
noiseless, well separated case (top) and an interference-rich, noisy
variant (bottom).
\label{fig:synreaction_datavis}}
\end{figure}

The resulting data matrix is obtained as the product of the kinetic-
and spectra matrix, viz. $M = WH$.  A visualization of the data matrix
$M$ is given in Figure~\ref{fig:synreaction_datavis} (top).

\subsection{Recovery in the noiseless case}
\label{sec:noiseless}

Given the data corresponding to the artificial reaction scheme
described in Section~\ref{sec:synthetic_reaction}, our goal is now to
recover both the component spectra of the five species and the
reaction kinetics \emph{only} from these data, i.e. to recover the
matrices $W$ and $H$ using only the data in $M$ without any further
information.  By construction of the data, we know that $M$ is
separable, and we will now apply the methods described in
Section~\ref{sec:sepnmf}.

\begin{figure}[p]
    \begin{center}
        \includegraphics[height=.45\textwidth]{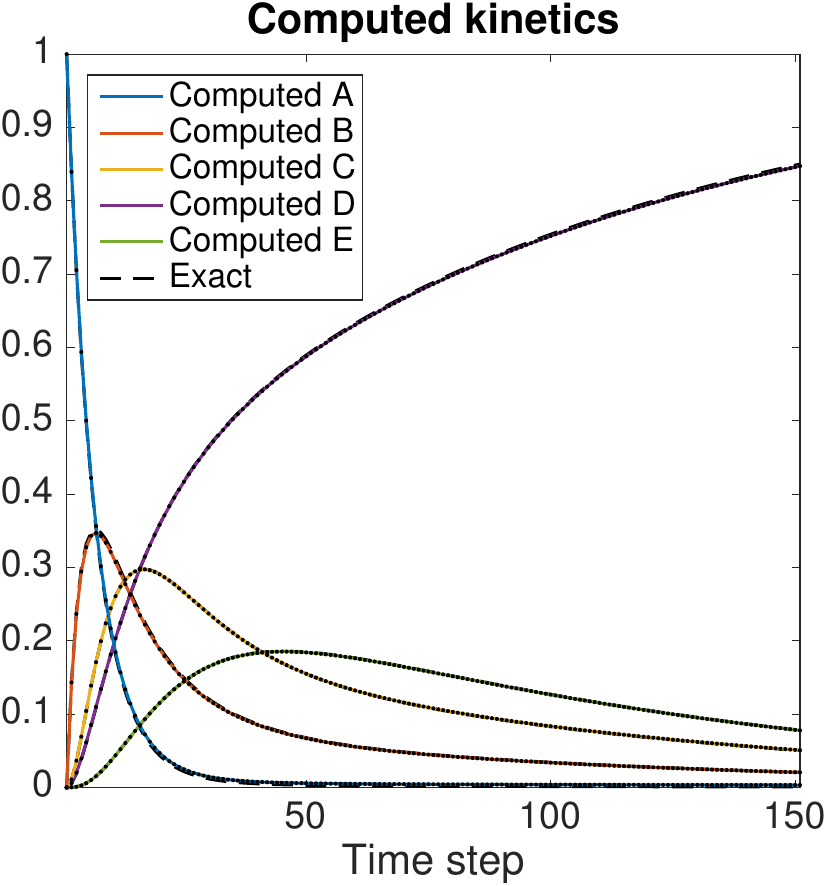}
        \includegraphics[height=.45\textwidth]{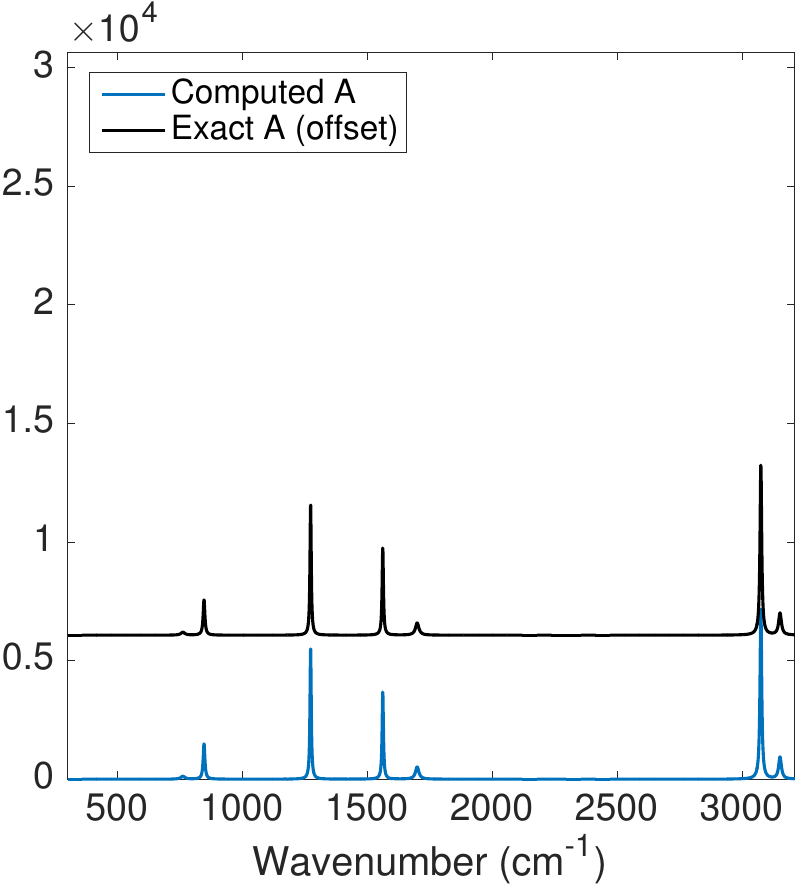}
        ~
        \includegraphics[height=.45\textwidth]{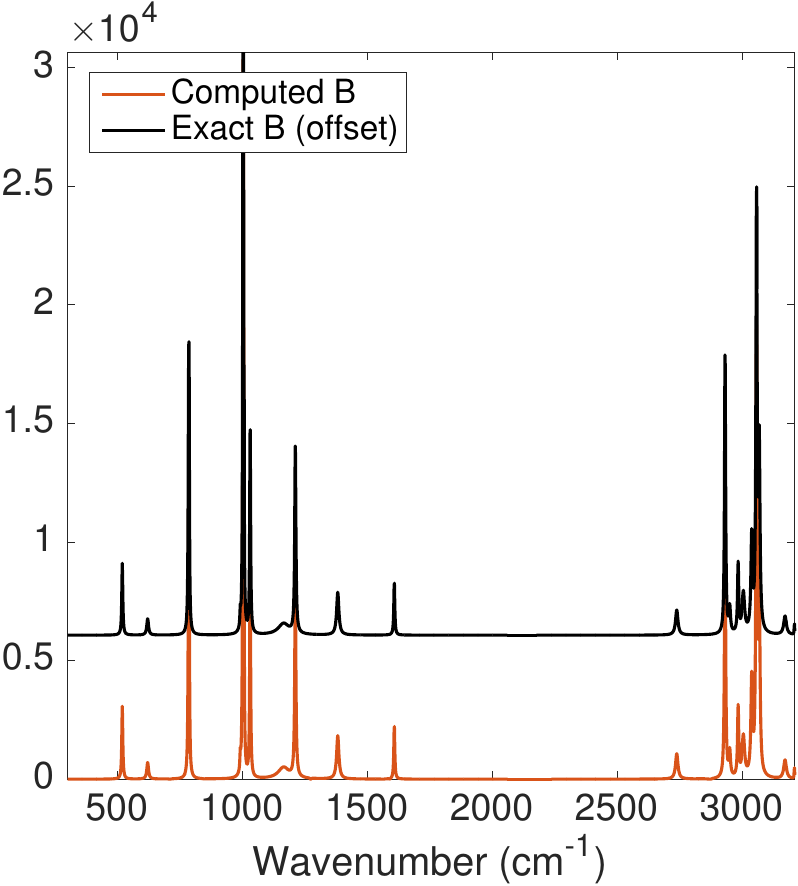}
        \includegraphics[height=.45\textwidth]{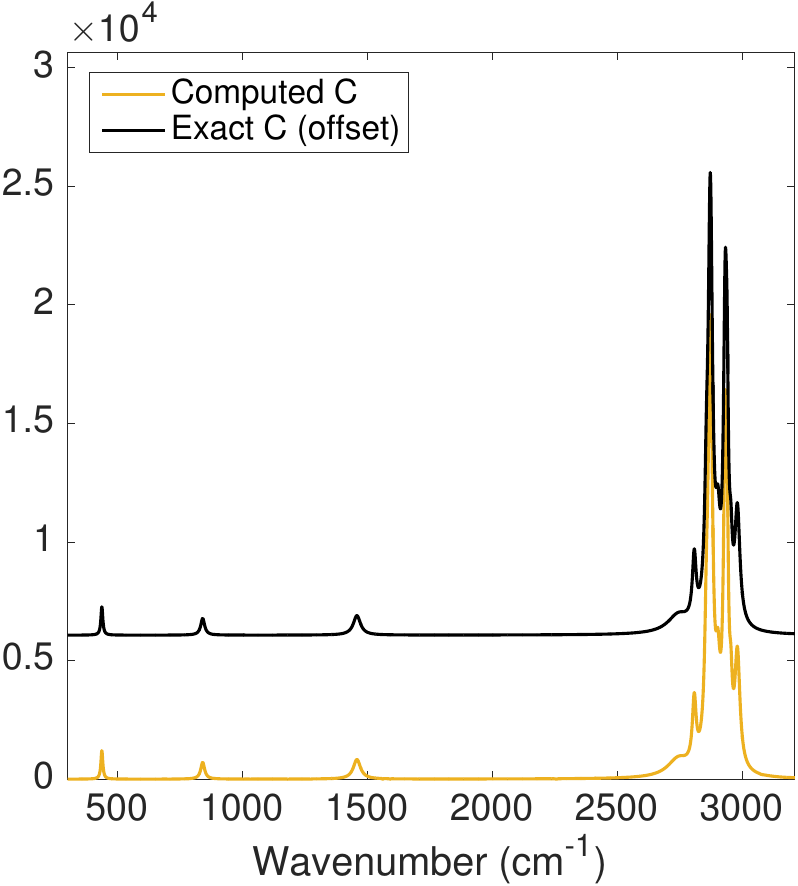}\\
        ~
        \includegraphics[height=.45\textwidth]{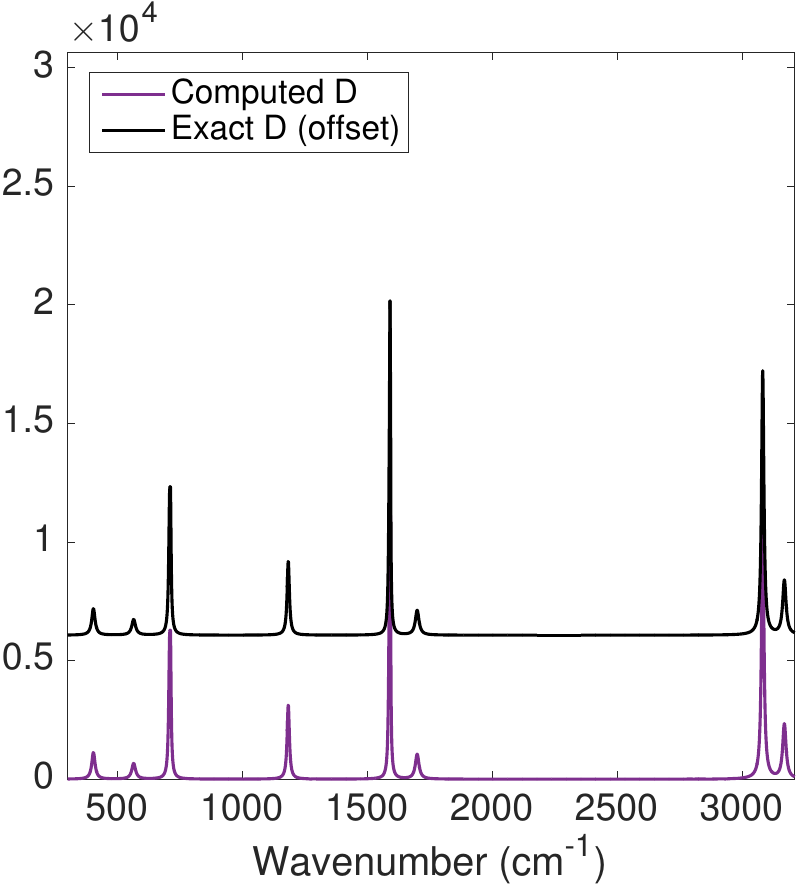}
        \includegraphics[height=.45\textwidth]{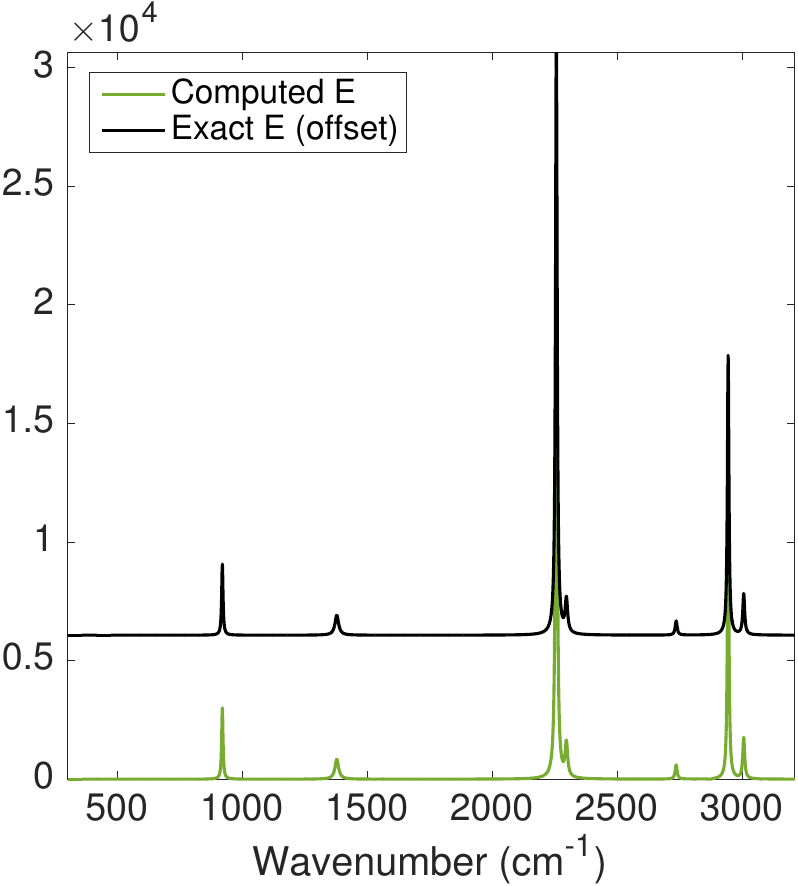}
    \end{center}
    \caption{Recovered reaction kinetics (top left) and spectral
fingerprints for the noiseless Raman measurements (see
Section~\ref{sec:noiseless}).  The computed and exact solutions are
visually almost indistinguishable.\label{fig:synreaction_n0}}
\end{figure}

The results are displayed in Figure~\ref{fig:synreaction_n0},
demonstrating that both the reaction kinetics $H$ and all the
component spectra are recovered to such a high level of accuracy, that
they can hardly be distinguished visually from the original data.

Nevertheless, the computed factors are not identical to the true
factors $W$ and $H$.  For the relative errors of $\tilde{W}$ and
$\tilde{H}$ we find
\begin{equation*}
    \frac{\norm{H - \tilde{H}}_F}{\norm{H}_F} = 0.0076, \quad
    \frac{\norm{W - \tilde{W}}_F}{\norm{W}_F} = 0.0051.
\end{equation*}
Hence the relative error for both the kinetics and the species is less
than 1\%.

Finally, we will recover the reaction coefficient matrix $K$ from the
computed kinetics matrix $\tilde{H}$.  Applying the methodology
described in Section~\ref{sec:comp_aspects}, we obtain
\begin{equation*}
\tilde{K} =
\begin{bmatrix}
   -0.5390 &  0.0349 & -0.0082 &  0.0003 &  0.0029\\
    0.5272 & -0.6580 &  0.2455 & -0.0005 &  0.0009\\
    0.0020 &  0.4295 & -0.3577 &  0.0001 &  0.0995\\
    0.0095 &  0.1941 &  0.0100 &  0.0001 & -0.0028\\
    0.0002 & -0.0006 &  0.1104 &  0.0000 & -0.1005
\end{bmatrix}.
\end{equation*}
Compared with the true reaction coefficients $K$ in
\eqref{eqn:k_matrix}, we find that the largest error made in
estimating any of the reaction coefficients from the computed kinetics
$\tilde{K}$ is $0.0159$, related to the coefficient $k_{12}$.

In summary, we find that all the data that constitute the reaction
network described in Section~\ref{sec:synthetic_reaction} have been
recovered quite accurately.

\subsection{Effect of increased interference}
\label{sec:move_bands}

In the model used in the previous section, the bands of the species
involved were (visually) well separated from each other.  In the other
extreme, i.e., if all bands of a species interfere with those of other
species, our approach will not be applicable for the recovery of the
reaction kinetics and the species fingerprints, as the factors are no
longer close to a separable factorization (the term $N_I$
in~\eqref{eqn:near_sep_interfere} becomes too large).  Thus we now
consider the case of ``modest'' interference.

\begin{figure}[t]
\begin{center}
\includegraphics[width=.47\textwidth]{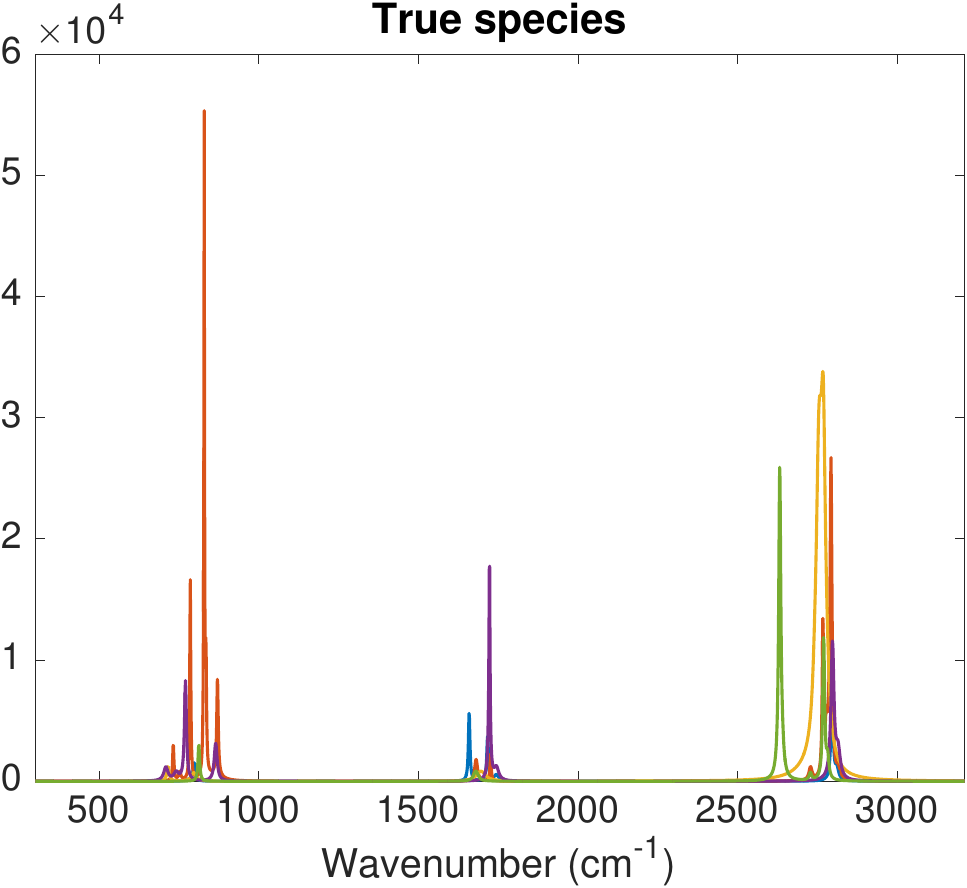}
\hfill \includegraphics[width=.49\textwidth]{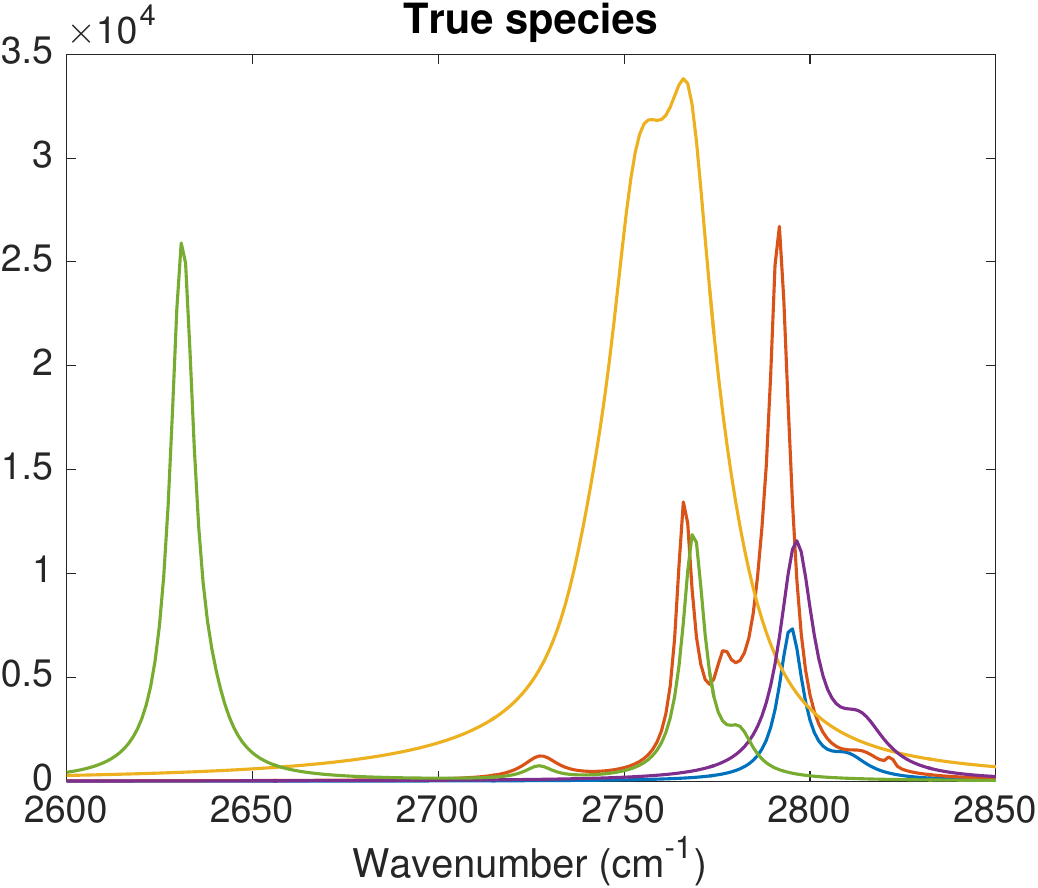}
\caption{{\itshape Left:} Species fingerprints where the Lorentz bands
have been moved closer to each other (compare with
Figure~\ref{fig:synreaction}). {\itshape Right:} Expanded view of the
high-frequency region. \label{fig:more_interference}}
\end{center}
\end{figure}

\begin{figure}[p]
    \begin{center}
        \includegraphics[height=.45\textwidth]{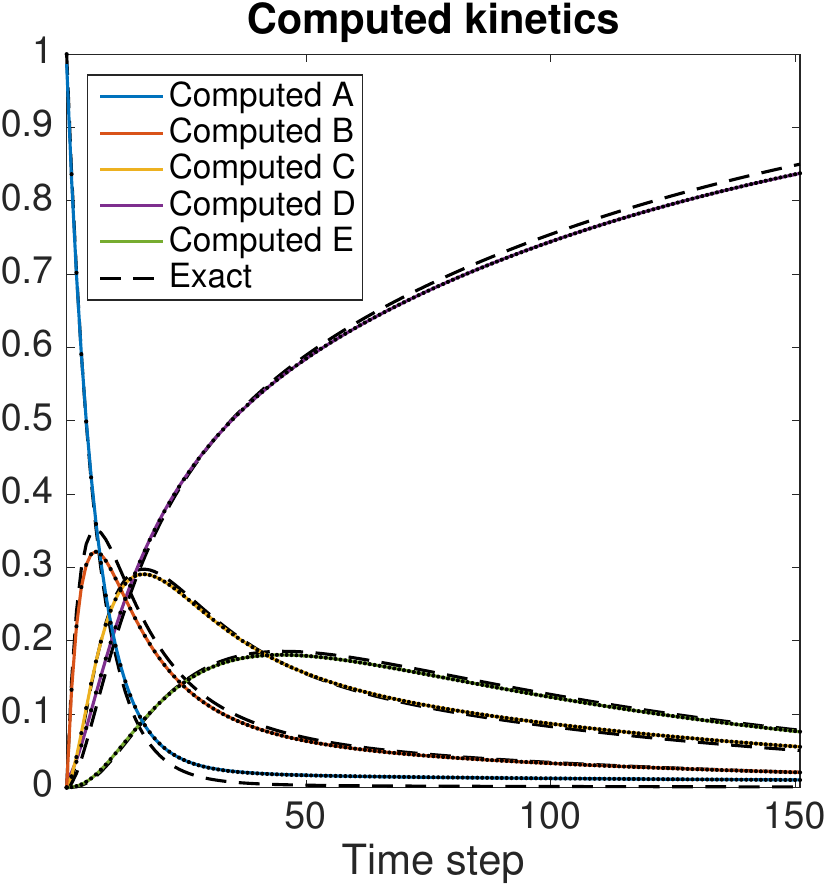}
        \includegraphics[height=.45\textwidth]{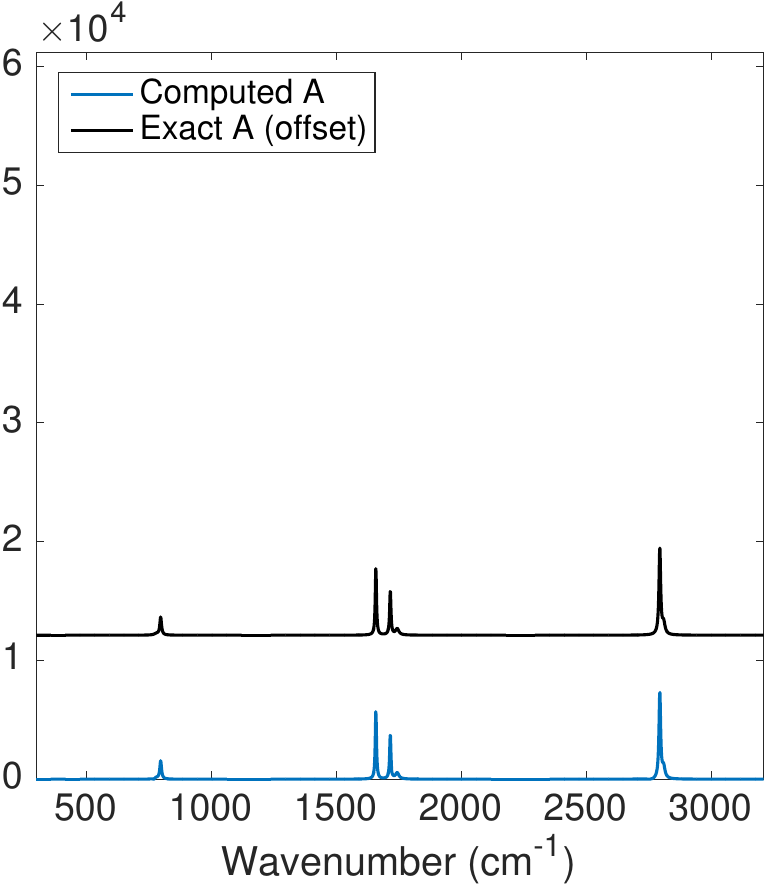}
        ~
        \includegraphics[height=.45\textwidth]{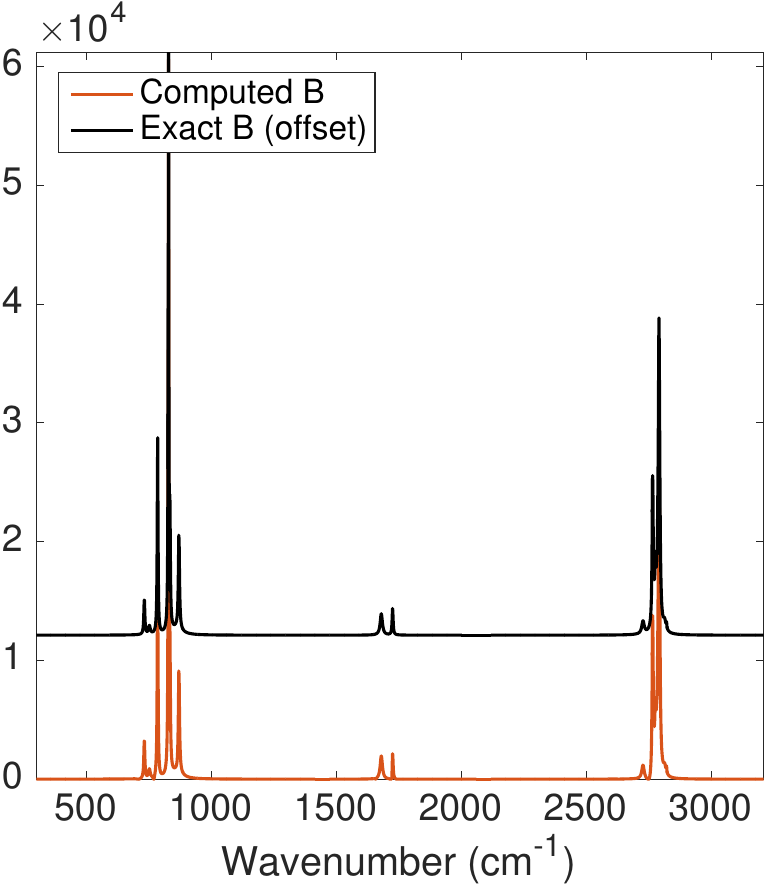}
        \includegraphics[height=.45\textwidth]{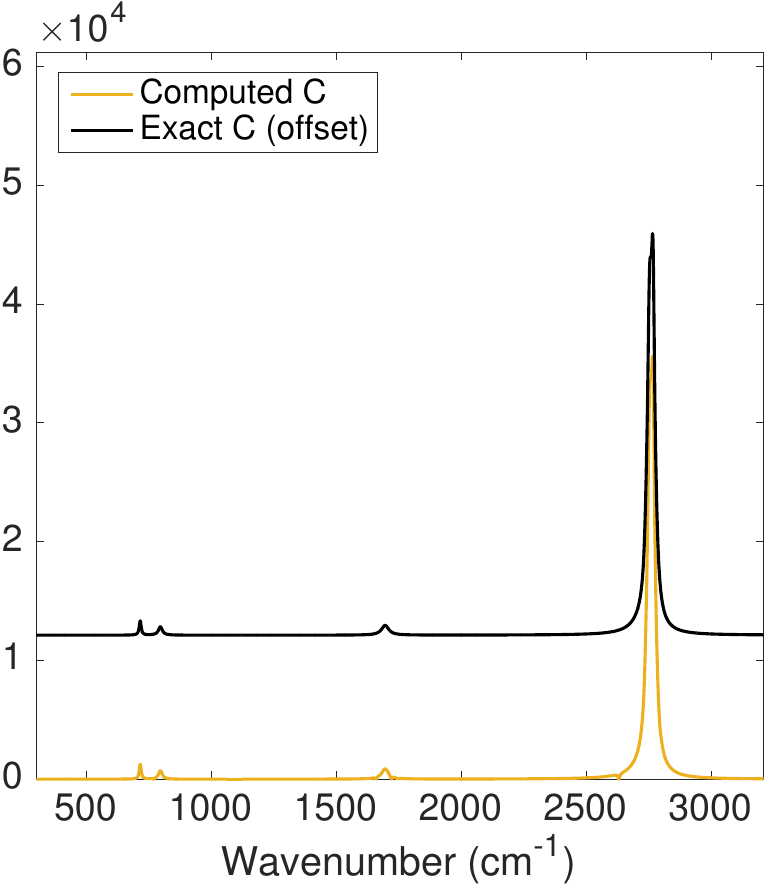}\\
        ~
        \includegraphics[height=.45\textwidth]{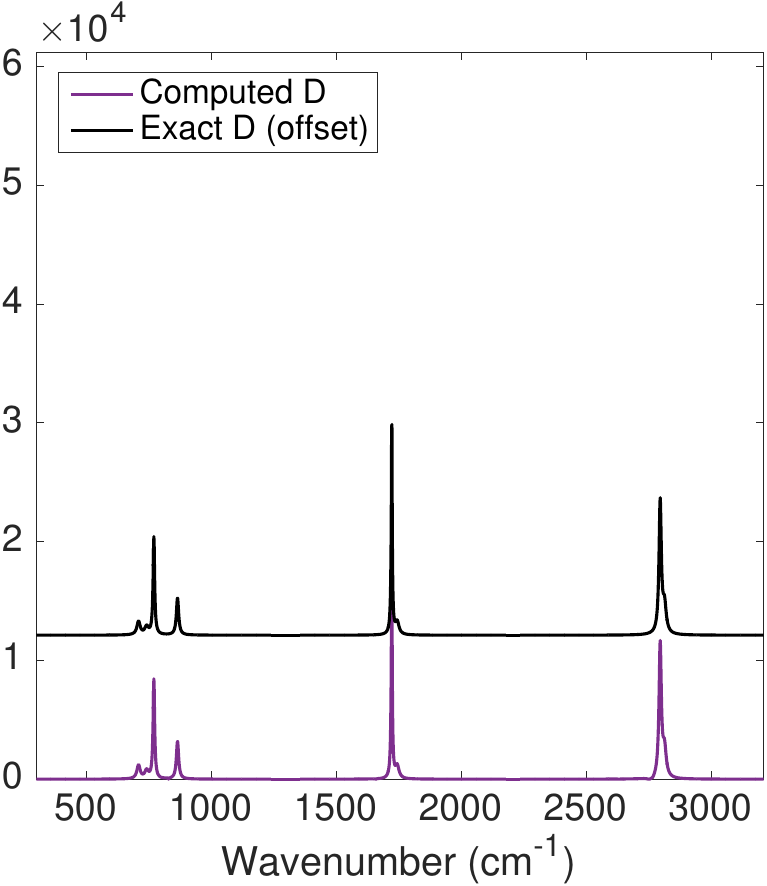}
        \includegraphics[height=.45\textwidth]{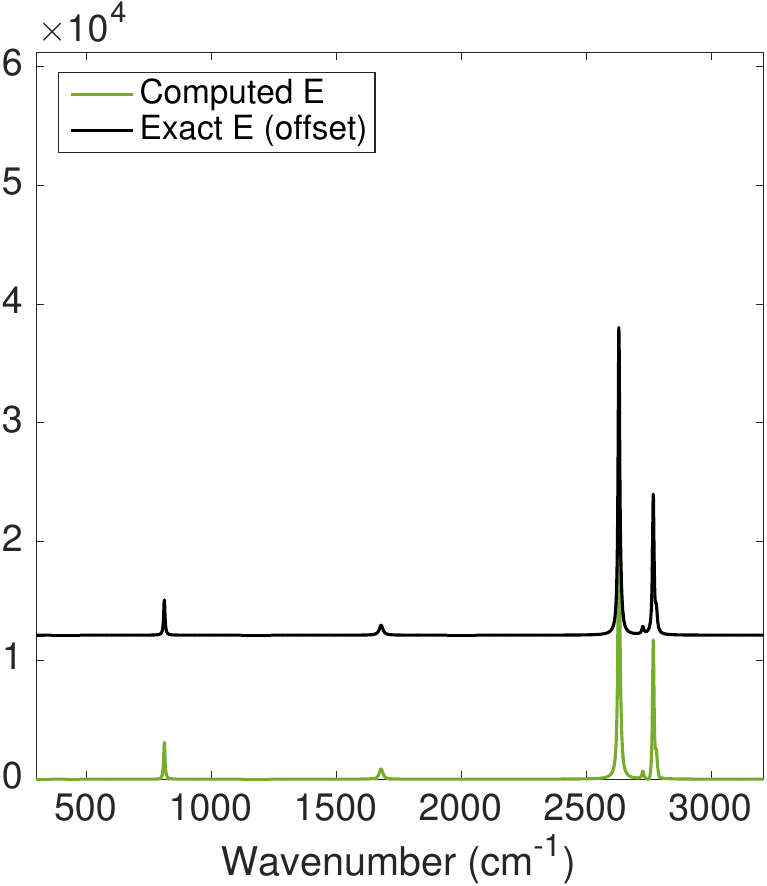}
    \end{center}
    \caption{Recovered reaction kinetics (top left) and spectral
fingerprints for the case of high interference (see
Section~\ref{sec:move_bands}).
The computed solutions are
still very good approximations to the true
solution.\label{fig:synreaction_interfere}}
\end{figure}

We enforce an increased level of interference among the species by
moving all the base points $x_0$ in all species towards three focal
points (see Figure~\ref{fig:more_interference}).
The result of Algorithm 1 being applied to these interference-rich
data is shown in Figure~\ref{fig:synreaction_interfere}.  Because of
the increased interference, the computed kinetics deviate slightly
from the true ones, but all the species bands have been identified
quite satisfactorily.  For the relative errors of $\tilde{W}$ and
$\tilde{H}$ we find
\begin{equation*}
    \frac{\norm{H - \tilde{H}}_F}{\norm{H}_F} = 0.063, \quad
    \frac{\norm{W - \tilde{W}}_F}{\norm{W}_F} = 0.026,
\end{equation*}
so the relative error for both the kinetics and the species spectra
is not more than $7\%$ and $3\%$, respectively.

If, however, the bands in the original species spectra are
moved even closer to each other, our method will eventually fail
to produce a qualitatively good solution.

\subsection{Recovery under the influence of measurement noise}
\label{sec:noisebench}

The data used in the previous section are highly idealized in the
sense that they are free of measurement noise.  In any practical
setting, the experimentally acquired Raman measurements will be
contaminated with noise from different sources, such as signal shot
noise or background noise (e.g. fluorescence).

We will now simulate the effect of noise added to the measurement
matrix $M$ in the sense of~\eqref{eqn:add_noise}.  We assume that all
the noise from different sources taken together resemble additive
Gaussian white noise.  Overall we disturb the data matrix
$M$ according to
\begin{equation*}
    \tilde{M} = M + \noiselevel\, \abs{N},
\end{equation*}
where the entries of $N \in \R^{m,n}$ are drawn from the normal
distribution $\mathcal{N}(0,1)$ and $\noiselevel = 0.4$ is the
relative noise level.  Further we will assume for our noise model,
that any constant background has already been removed from the
measurement data.  In Figure~\ref{fig:synreaction_datavis} (bottom),
we visualize the resulting noisy data $\tilde{M}$.

\begin{figure}[p]
    \begin{center}
        \includegraphics[height=.45\textwidth]{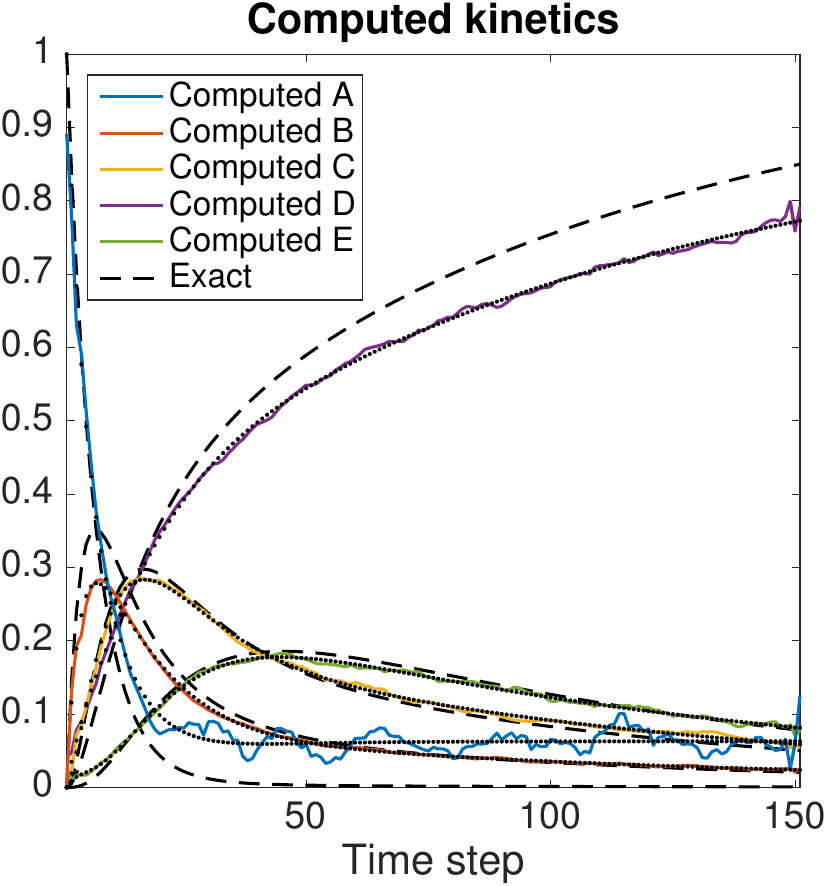}
        \includegraphics[height=.45\textwidth]{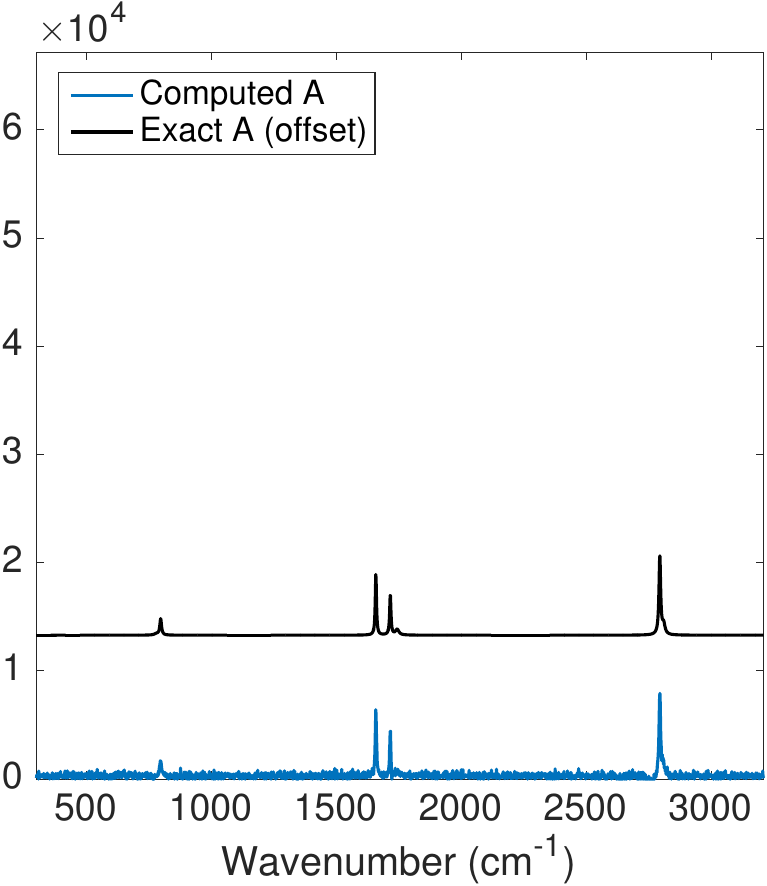}
        ~
        \includegraphics[height=.45\textwidth]{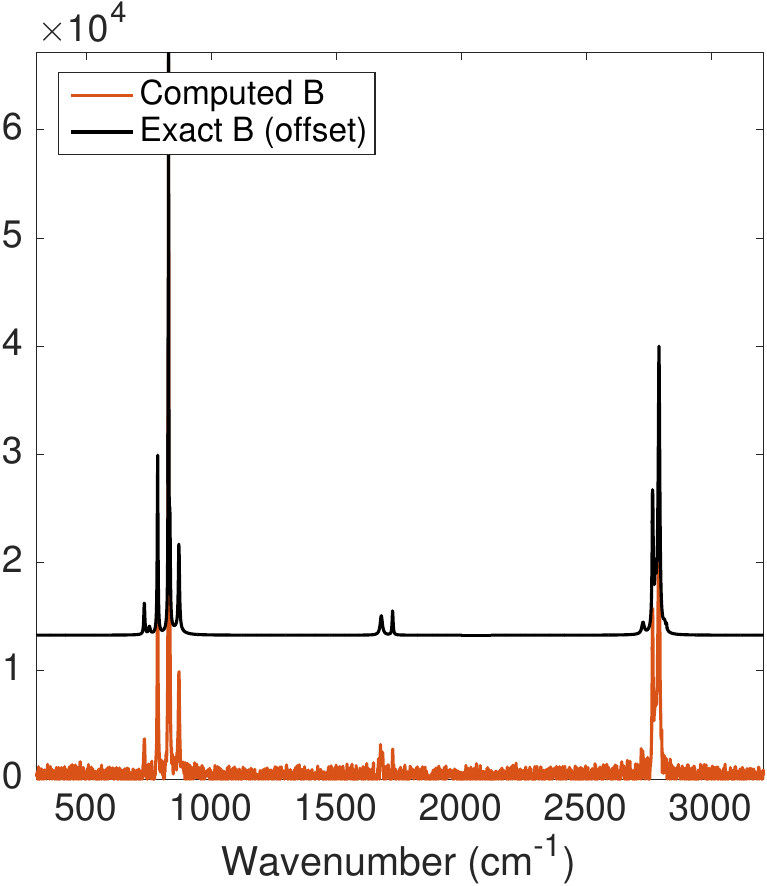}
        \includegraphics[height=.45\textwidth]{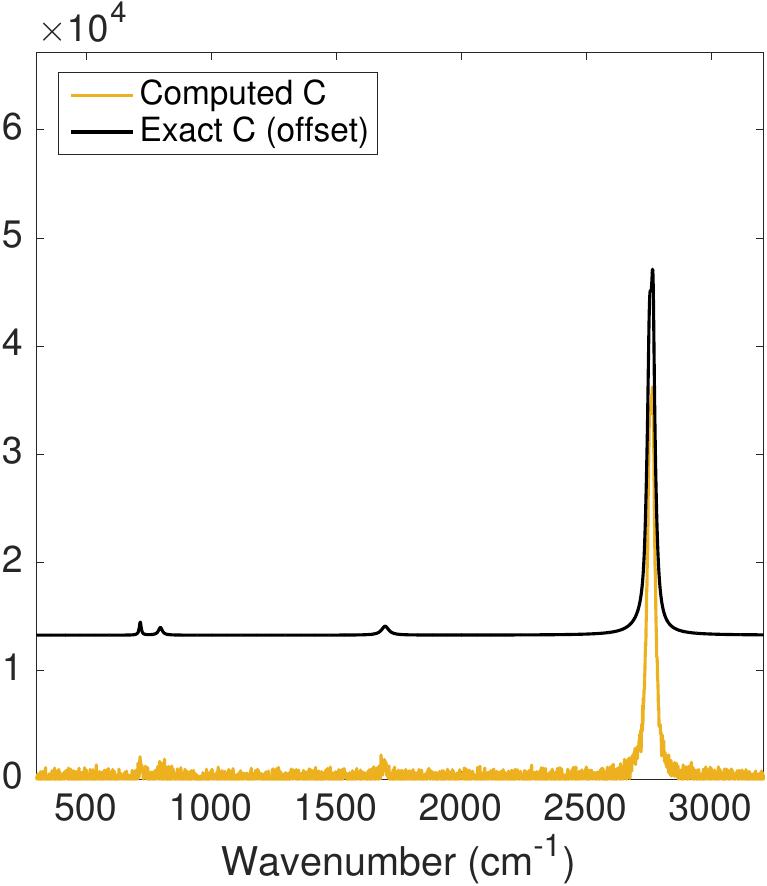}\\
        ~
        \includegraphics[height=.45\textwidth]{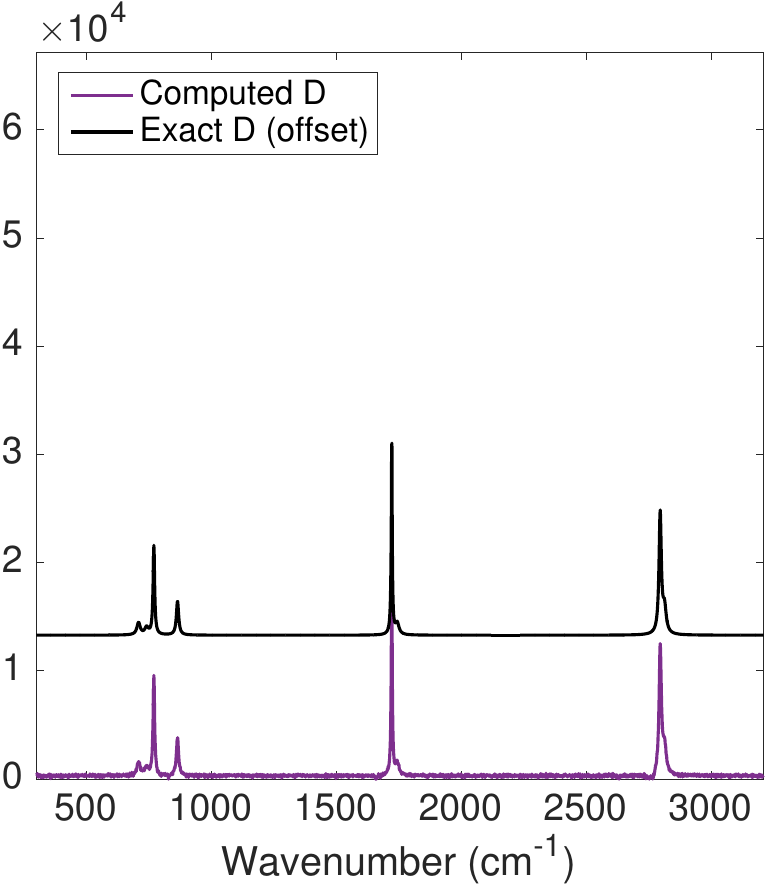}
        \includegraphics[height=.45\textwidth]{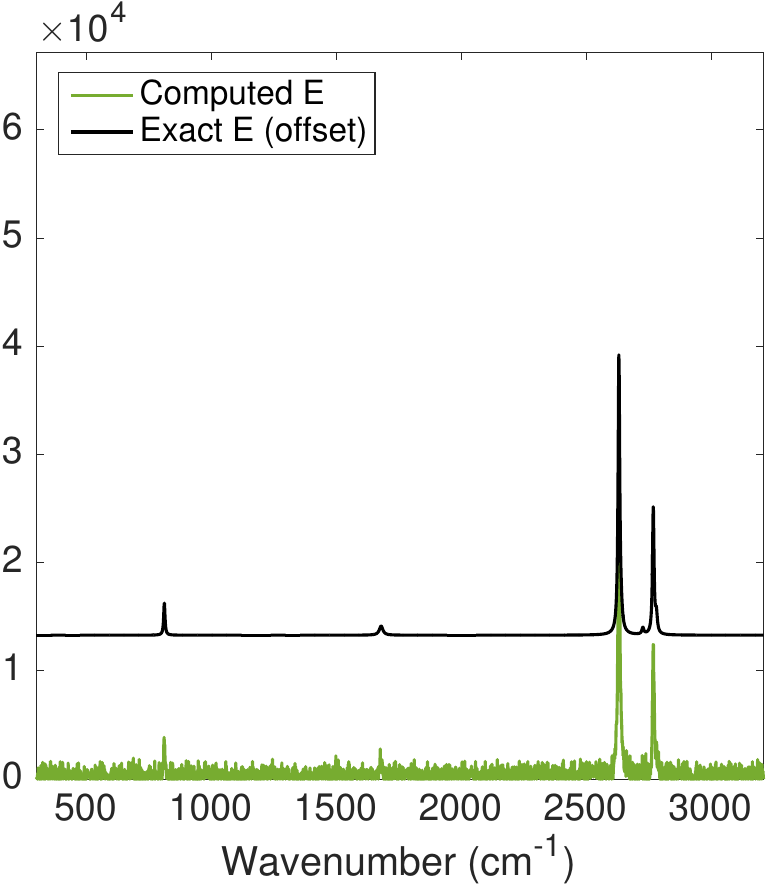}
    \end{center}
    \caption{Recovered reaction kinetics (top left) and spectral
fingerprints for the noisy Raman measurements (see
Section~\ref{sec:noisebench}).  The dotted lines in the reaction
kinetic show the kinetic corresponding to the extracted reaction
coefficients, see Step $6$ in Algorithm 1.\label{fig:synreaction_noisy}}
\end{figure}

The component spectra and kinetics were recovered as described in
Section~\ref{sec:overall_algo}.  Figure~\ref{fig:synreaction_noisy}
shows the result of our algorithm applied to the noisy,
interference-rich data.  Because of the large distance of $W,H$ to
being truly separable, the computed kinetics displays some deviations
from the true data, but still provides a satisfactory description.
Also the computed component spectra show a reasonable agreement with
the true spectra.

\subsection{Determination of the number of species}

\label{sec:err_comp}

In the numerical experiments described above we have always assumed
that the number of species is known (here $r=5$). Of course, in
practical applications the correct determination of the value of $r$
is one of the greatest challenges. One approach for solving this
problem is to use Algorithm~1 for different values of $r$, which
yields approximate species $W_r$ and kinetics $H_r$, and then compute
the (relative) data error $\norm{M-W_r H_r}_F / \norm{M}_F$.

In the following table we show the data errors for
$r=1,2,\dots,7$ and the three experimental setups from
Sections~4.2--4.4: noiseless case (first row of the table), increased
interference (second), increased interference and measurement noise
(third). In each case a significant drop in the data error occurs at
the correct number of species, while no significant further reduction of the
error is achieved when increasing the number of (suspected) species
even further.

\begin{table}[h]
\begin{center}
\begin{small}
\begin{tabular}{c|ccccccc}
\label{table:data_error}
{$r$}  & 1 & 2 & 3 & 4 & 5 & 6 & 7\\
\hline
{Sec.~4.2} &
0.97294 &
0.77203 &
0.18479 &
0.06055 &
0.00003 &
0.00003 &
0.00003 \\
{Sec.~4.3} &
0.94908 &
0.68692 &
0.20209 &
0.04698 &
0.00029 &
0.00029 &
0.00029 \\
{Sec.~4.4} &
0.83106 &
0.35383 &
0.25894 &
0.12683 &
0.08169 &
0.08107 &
0.08066 
\end{tabular}
\end{small}
\caption{Relative data fit errors when applying Algorithm~1 to the
setups from Sections~\ref{sec:noiseless}--\ref{sec:noisebench} and
different values of $r$.}
\end{center}
\end{table}

The results indicate that using the NMF for determining the number of
species is an alternative to existing techniques in this context such
as singular values of the data matrix.  An extensive survey of the
latter technique is given in~\cite{MalBook02}.

\section{Concluding remarks}
\label{sec:concl}

We have presented an algorithm for the recovery of component spectra
and reaction kinetics from data obtained through time resolved Raman
spectroscopy. The key tool we used in our approach is the
so-called separability condition in the non-negative matrix factorization
(NMF). In terms of the component spectra this condition is
approximately satisfied, if each of the species has at least one
band in its spectrum that does not interfere too much with the bands
of other species. Thus, it is conceptually similar to the classical
Rayleigh criterion limit. Our approach combines a standard algorithm
for separable NMF (``SNPA'') with a few
pre-processing steps for the measurement data. A number of other
separable NMF algorithms exist, but we did not yet pursue a detailed
comparison among them for the present application. In our numerical
study we have demonstrated that the component spectra and reaction
kinetics can be recovered with a reasonable quality under modest
measurement noise and interference among the component spectra.
Whereas in terms of the kinetics the approach currently restricted to
a network of first-order or pseudo-first-order reactions, it may
equally well applied to other spectroscopic techniques such as IR or
NMR spectroscopies.

\paragraph*{Acknowledgements}  This research was supported by the DFG
cluster of excellence ``UniCat''.

\bibliographystyle{applspec}
\bibliography{using_sepnmf}

\begin{thebibliography}{10}
\newcommand{\enquote}[1]{``#1''}

\bibitem{AndoHamaguchi2013}
M.~Ando and H.-o. Hamaguchi, \enquote{Molecular Component Distribution Imaging
  of Living Cells by Multivariate Curve Resolution Analysis of Space-Resolved
  Raman Spectra}. J. Biomed. Opt. 2013. 19(1):011016.

\bibitem{AroraEtAl:2012}
S.~Arora, R.~Ge, R.~Kannan, and A.~Moitra, \enquote{Computing a Nonnegative
  Matrix Factorization---Provably}. In S{TOC}'12---{P}roceedings of the 2012
  {ACM} {S}ymposium on {T}heory of {C}omputing, pages 145--161, ACM, New York,
  2012.

\bibitem{BalakrishnanEtAl2008}
G.~Balakrishnan, C.~L. Weeks, M.~Ibrahim, A.~V. Soldatova, and T.~G. Spiro,
  \enquote{{Protein Dynamics from Time Resolved UV Raman Spectroscopy.}} Curr.
  Opin. Struct. Biol. 2008. 18(5):623--629.

\bibitem{BeljebbarEtAl2009}
A.~Beljebbar, O.~Bouch{\'e}, M.~D. Di{\'e}bold, P.~J. Guillou, J.~P. Palot,
  D.~Eudes, and M.~Manfait, \enquote{{Identification of Raman Spectroscopic
  Markers for the Characterization of Normal and Adenocarcinomatous Colonic
  Tissues.}} Crit. Rev. Oncol. Hematol. 2009. 72(3):255--264.

\bibitem{Bioucas-DiasEtAl2012}
J.~M. Bioucas-Dias, A.~Plaza, N.~Dobigeon, M.~Parente, Q.~Du, P.~Gader, and
  J.~Chanussot, \enquote{{Hyperspectral Unmixing Overview: Geometrical,
  Statistical, and Sparse Regression-Based Approaches}}. IEEE J. Sel. Top.
  Appl. Earth Observations Remote Sensing 2012. 5(2):354--379.

\bibitem{ChanEtAl2008}
T.-H. Chan, W.-K. Ma, C.-Y. Chi, and Y.~Wang, \enquote{A Convex Analysis
  Framework for Blind Separation of Non-Negative Sources}. IEEE Trans. Signal
  Process. 2008. 56(10, part 2):5120--5134.

\bibitem{DonohoStodden:2003}
D.~Donoho and V.~Stodden, \enquote{When Does Non-Negative Matrix Factorization
  give a Correct Decomposition into Parts?} In S.~Thrun, L.~Saul, and
  B.~{Sch\"{o}lkopf}, editors, Advances in Neural Information Processing
  Systems 16, MIT Press, Cambridge, MA, 2004.

\bibitem{DoepnerEtAl1996}
S.~Döpner, P.~Hildebrandt, A.~G. Mauk, H.~Lenk, and W.~Stempfle,
  \enquote{Analysis of Vibrational Spectra of Multicomponent Systems.
  {A}pplication to p{H}-dependent Resonance {R}aman Spectra of ferricytochrome
  c}. Spectrochim. Acta, Part A 1996. 52(5):573 -- 584.

\bibitem{Gillis:2013}
N.~Gillis, \enquote{Robustness {A}nalysis of {H}ottopixx, a {L}inear
  {P}rogramming {M}odel for {F}actoring {N}onnegative {M}atrices}. SIAM J.
  Matrix Anal. Appl. 2013. 34(3):1189--1212.

\bibitem{Gillis2014}
N.~Gillis, \enquote{Successive Nonnegative Projection Algorithm for Robust
  Nonnegative Blind Source Separation}. SIAM J. Imaging Sci. 2014.
  7(2):1420--1450.

\bibitem{Gillis:HowAndWhy}
N.~Gillis, \enquote{{The Why and How of Nonnegative Matrix Factorization}}. In
  Chapman {\&} Hall/CRC Machine Learning {\&} Pattern Recognition, pages
  257--291, Chapman and Hall/CRC, 2014.

\bibitem{GillisLuce:2013}
N.~Gillis and R.~Luce, \enquote{Robust Near-Separable Nonnegative Matrix
  Factorization using Linear Optimization}. J. Mach. Learn. Res. 2014.
  15:1249--1280.

\bibitem{GillisVavasis:2012}
N.~Gillis and S.~A. Vavasis, \enquote{Fast and Robust Recursive Algorithms for
  Separable Nonnegative Matrix Factorization}. IEEE Trans. Pattern Anal. Mach.
  Intell. 2014. 36(4):698--714.

\bibitem{HendlerShrager1994}
R.~W. Hendler and R.~I. Shrager, \enquote{{Deconvolutions Based on Singular
  Value Decomposition and the Pseudoinverse: a Guide for Beginners.}} J.
  Biochem. Biophys. Methods 1994. 28(1):1--33.

\bibitem{HenryHofrichter1992}
E.~Henry and J.~Hofrichter, \enquote{Singular Value Decomposition: Application
  to Analysis of Experimental Data}. In Numerical Computer Methods, volume 210
  of \emph{Methods in Enzymology}, pages 129 -- 192, Academic Press, 1992.

\bibitem{KumarSindhwaniPrabhanjan2013}
A.~Kumar, V.~Sindhwani, and P.~Kambadur, \enquote{Fast Conical Hull Algorithms
  for Near-Separable Non-Negative Matrix Factorization}. ICML 2013. (PART
  1):231--239.

\bibitem{LeeSeung:1999}
D.~D. Lee and H.~S. Seung, \enquote{Learning the Parts of Objects by
  Non-Negative Matrix Factorization}. Nature 1999. 401(6755):788--791.

\bibitem{MalBook02}
E.~R. Malinkowski, Factor Analysis in Chemistry. John Wiley \& Sons, Inc., New
  York, 3rd edition, 2002.

\bibitem{Sawall2010}
K.~Neymeyr, M.~Sawall, and D.~Hess, \enquote{Pure Component Spectral Recovery
  and Constrained Matrix Factorizations: Concepts and Applications}. J. Chemom.
  2010. 24(2):67--74.

\bibitem{OndriasSimpsonLarson1996}
M.~Ondrias, M.~Simpson, and R.~Larsen, \enquote{Time-Resolved Resonance Raman
  Spectroscopy}. In J.~Laserna, editor, Modern Techniques in Raman
  Spectroscopy, Wiley, 1996.

\bibitem{SahooUmapathyParker2011}
S.~K. Sahoo, S.~Umapathy, and A.~W. Parker, \enquote{{Time-Resolved Resonance
  Raman Spectroscopy: Exploring Reactive Intermediates}}. Appl. Spectrosc.
  2011. 65(10):1087--1115.

\bibitem{ShinzawaEtAl2009}
H.~Shinzawa, K.~Awa, W.~Kanematsu, and Y.~Ozaki, \enquote{{Multivariate Data
  Analysis for Raman Spectroscopic Imaging}}. J. Raman Spectrosc. 2009.
  40(12):1720--1725.

\bibitem{TalariEtAl2015}
A.~C.~S. Talari, C.~A. Evans, I.~Holen, R.~E. Coleman, and I.~U. Rehman,
  \enquote{{Raman Spectroscopic Analysis Differentiates Between Breast Cancer
  Cell Lines}}. J. Raman Spectrosc. 2015. 46(5):421--427.

\bibitem{Thomas:1974}
L.~Thomas, \enquote{Rank Factorization of Nonnegative Matrices ({A}.
  {B}erman)}. SIAM Rev. 1974. 16(3):393--394.

\bibitem{VanBenthemKeenan:2004}
M.~H. Van~Benthem and M.~R. Keenan, \enquote{Fast Algorithm for the Solution of
  Large-Scale Non-Negativity-Constrained Least Squares Problems}. J. Chemom.
  2004. 18(10):441--450.

\bibitem{Vavasis:2008}
S.~A. Vavasis, \enquote{On the Complexity of Nonnegative Matrix Factorization}.
  SIAM J. Optim. 2009. 20(3):1364--1377.

\bibitem{WeakleyEtAl2012}
A.~T. Weakley, {Warwick, P. C. Temple}, T.~E. Bitterwolf, and D.~E. Aston,
  \enquote{{Multivariate Analysis of Micro-Raman Spectra of Thermoplastic
  Polyurethane Blends Using Principal Component Analysis and Principal
  Component Regression}}. Appl. Spectrosc. 2012. 66(11):1269--1278.

\bibitem{ZhangEtAl2005}
A.~Zhang, W.~Zeng, T.~M. Niemczyk, M.~R. Keenan, and D.~M. Haaland,
  \enquote{{Multivariate analysis of Infrared Spectra for Monitoring and
  Understanding the Kinetics and Mechanisms of Adsorption Processes}}. Appl.
  Spectrosc. 2005. 59(1):47--55.

\end{thebibliography}

\end{document}